\documentclass[12pt,notitlepage,a4paper,reqno,ps]{amsart}
\usepackage{eurosym}
\usepackage{amssymb}
\usepackage{amstext}
\usepackage{color}

\newcommand{\medint}{-\kern  -,375cm\int}

\newenvironment{michelarev}{\color{red}}{\color{black}}
\newcommand{\bmicr}{\begin{michelarev}}
\newcommand{\emicr}{\end{michelarev}}
\theoremstyle{plain}
\newtheorem{theorem}{Theorem}

\theoremstyle{definition}

\theoremstyle{remark}
\newtheorem{remark}[theorem]{Remark}

\makeatother  \allowdisplaybreaks
\email{Corresponding author: michela.eleuteri@unimore.it (Michela Eleuteri)}

\begin{document}
\title[A Volterra's generalization applied to flat functions]{The Fundamental Theorem of Integral Calculus: \\ a Volterra's generalization applied to flat functions}
\author[C. Benassi, M. Eleuteri]{Carlo Benassi, Michela Eleuteri}
\address{Dipartimento di Scienze Fisiche, Informatiche e Matematiche,
via Campi 213/b, 41125 Mo-dena, Italy}

\thanks{
The authors are partially supported GNAMPA (Gruppo Nazionale per l'Analisi Matematica, la Probabilit\`a e le loro Applicazioni)
 of INdAM (Istituto Nazionale di Alta Matematica) and by the University of Modena and Reggio Emilia.}

\begin{abstract}
In a recent paper \cite{S18} a smooth function $f: [0,1] \rightarrow \mathbb{R}$ with all derivatives vanishing at 0 has been considered and a global condition, showing that $f$ is indeed identically 0, has been presented. The purpose of this note is to replace the classical Fundamental Theorem of Calculus for the Riemann integral, as it has been used in \cite{S18}, with a weaker form going back to Volterra \cite{V}, which is little known. Therefore the proof we propose in this paper turns to be important also from the teaching point of view, as long as in literature there are very few examples in which explicitly the lower integral and the upper integral of a function appear (usually the assumption that the function is Riemann-integrable is required).
\end{abstract}

\maketitle

\begin{center}
\fbox{\today}
\end{center}

It is well known that one of the main reasons why Lebesgue integral replaced the Riemann one in the applications is the fact that in the theory of Lebesgue integration the Fundamental Theorem of Calculus holds under weaker assumptions; in particular all absolutely continuous functions can be reconstructed starting from their derivative, by means of the Lebesgue integral.
\\
In the theory of Riemann integration the most general version of the  Fundamental Theorem of Calculus is the following:
\begin{theorem}
\label{R-general}
Let $f: [a,b] \rightarrow \mathbb{R}$ be a real valued differentiable function such that $f'$ is integrable; then, for all $x \in [a,b]$ it holds
\[
f(x) = f(a) + \int_a^x f'(t) \, dt.
\]
\end{theorem}
When $f'$ is not integrable, the Riemann integration theory allows anyway, sometimes, to get information on $f,$ without having to resort to more sophisticated notion of integration, such as for instance the Lebesgue or the Henstock one.
\\
The aim of this paper is just that of showing a situation of this kind, which we believe it may be of interest because the Riemann integral has retained a central role from the didactic point of view.
\\
In particular, the main result of this paper is Theorem \ref{main-result}, which generalizes a result by G. Stoica \cite{S18}. In order to prove Theorem \ref{main-result}, we will use a very nice theorem by V. Volterra \cite{V}, see also \cite{PS} pag. 48, that is (undeservedly) little known. We notice that Theorem \ref{R-general} is a corollary of Theorem \ref{Volterra}.
\\
\begin{theorem}
\label{Volterra}
Let $h: [a,b] \rightarrow \mathbb{R}$ be a bounded function. Suppose that there exists $H: [a,b] \rightarrow \mathbb{R}$ such that $H'(x) = h(x)$ for all $x \in [a,b]$. Then
\[
\int_{{\!\!\!\!\!\!\!\!\!*} \, {\,\,\,\,a}}^b h(x) \, dx \le H(b) - H(a) \le  \int_a^{{\!\!\!\!\!\!\!*} \, {\,\,\,\,b}} h(x) \, dx.
\]
\end{theorem}
Here we denote with $\,\,\,\, \int_{{\!\!\!\!\!\!\!\!\!*} \, {\,\,\,\,a}}^b h(x) \, dx \,\,$ the {\it lower integral} and with $\,\,\,\, \int_a^{{\!\!\!\!\!\!\!*} \, {\,\,\,\,b}} h(x) \, dx \,\,$ the {\it upper integral} in the Riemann formulation. 
\\
\\
A {\it flat} function $f: [0,1] \rightarrow \mathbb{R}$ is a smooth function with all derivatives vanishing at zero. At the beginning of the last century Denjoy and Carleman (see for instance \cite{C68} for more details) considered the class of  {\it quasi-analytic} functions, being flat at 0, infinitely differentiable on $[0,1]$ and satisfying suitable inequalities involving all the derivatives, for instance
\begin{equation}
\label{carleman}
|f^{(n)}(x)| \le \, m_n \, a^{n+1},
\end{equation}
for some $a > 0$, $m_n \ge 0$ and for every $n \in \mathbb{N}$ and $x \in [0,1]$. They were able to characterize those sequences $\{m_n\}_{n \in \mathbb{N}}$ for which $f$ is identically zero. These conditions remained of a fundamental importance also in the works by Bernstein and Mandelbrojt (see for instance \cite{B64}, \cite{M52}) in the second half of the century.

In a very recent paper \cite{S18}, condition \eqref{carleman} has been replaced by a global simpler inequality, namely
\[
|x \, f'(x)| \le \, C \, |f(x)|
\]
for some $C > 0$ and every $x \in [0,1]$. This condition involves a $\mathcal{C}^1-$function and its first derivative and, together with the assumption that $f^{(n)}(0) = 0$ for every $n \in \mathbb{N}$, can be used to conclude that $f$ is identically 0 in $[0,1]$. The assumption that $f$ is a $\mathcal{C}^1$ function turns to be essential for the author in order to apply the Fundamental Theorem of Calculus for the Riemann integral. 
\\
In the next theorem, which constitutes our main result, we will show that actually it is not necessary to require that $f'$ is a Riemann-integrable function.
It can be stated as follows:
\begin{theorem}
\label{main-result}
Let $f: [0,1] \rightarrow \mathbb{R}$ be a real valued differentiable function such that
\begin{equation}
\label{zero}
 f(0) = 0 
\end{equation}
and
\begin{equation}
\label{uno}
\forall n > 1  \,\,\, \exists \delta_n > 0: \,\,\, |f(x)| < x^n \qquad \forall x \in ]0, \delta_n[.
\end{equation} 
Assume moreover that
\begin{equation}
\label{stoica}
\exists C > 0: \,\,\, |x f'(x)| \le \, C |f(x)| \qquad \forall x \in [0,1]. 
\end{equation}
Then $f(x) = 0$ for every $x \in [0,1]$.
\end{theorem}
In order to prove this result we will apply Theorem \ref{Volterra} with $h = |f'(x)|$. We will prove therefore that \eqref{zero}--\eqref{stoica} entails that $h$ is bounded so we are not assuming any integrability on our function $h$. Condition \eqref{zero} actually turns out to be equivalent to $f^{(n)}(0) = 0$ if all these derivatives exist.
\\
\begin{remark}
As we already mentioned before, what we believe it is interesting to remark is that in the proof of Theorem \ref{main-result} only the Riemann integral calculus has been employed. However, the proof could have been simplified if we would have used a more powerful tool. For instance, by means of the Henstock-Kurzweill integral \cite{G94}, Volterra's theorem would have been replaced by the Fundamental Theorem of Calculus (if we interpret the integrals in the sense of Henstock-Kurzweil, given a real valued differentiable function $h$, one always has that $\displaystyle \int_a^b h'(s) \, ds = h(b) - h(a)$ \Bigg).
\end{remark}

\begin{proof} First of all we observe that $|f'|$ is bounded. 
Indeed, condition \eqref{zero} and \eqref{stoica} imply that $\frac{f(x)}{x}$ is bounded (it is a continuous function which has limit in 0) therefore also $f'$ is also bounded. 
\\
Let us set 
\[
g(x) := |f(x)|
\]
Obviously $g$ is differentiable in all points $x$ for which $f(x) \neq 0$. On the other hand, if $x$ is such that $f(x) = 0$, then also $g$ turns to be differentiable because in these points $f'(x) = 0$ (this fact can be deduced from \eqref{uno}, if $x = 0$ and from \eqref{stoica} if $x \in ]0, 1]$).
Therefore $g$ is differentiable in the whole interval $[0,1]$ and in particular we have 
\[
|g'(x)| = |f'(x)| \qquad \forall x \in [0,1];
\]
thus also $|g'(x)|$ is bounded.

At this point, by Theorem \ref{Volterra} we show that $f \equiv 0$ on the whole interval $[0,1]$.
\\ 
We would like to show that 
\[
\forall n > C,  \qquad |f(x)| < x^n \qquad \forall x \in ]0,1]. 
\]
Notice that condition \eqref{uno} holds only in a small interval $]0, \delta_n[$. Notice that the constant $C$ is the one appearing in \eqref{stoica}. 
\\
Suppose by contradiction that this does not hold, i.e. suppose that there exists $n > C$ such that $f(x) \ge x^n$ for all $x \in ]0,1]$. By \eqref{uno} this means that there exists $n > C$ such that the set
\[
E_n := \{x > 0: \, |f(x)| = x^n\}
\]
is non-empty. From \eqref{uno} we deduce that $E_n \subset [\delta_n, 1]$; moreover $E_n$ is closed as long as $f$ is a continuous function, therefore it has a minimum point 
\[
\bar{x}_n := \min E_n > 0.
\]
Then, in $]0, x_n[,$ we have $|f(x)| < x^n$ and so, by Theorem \ref{Volterra} applied to the function $g$ 
\begin{eqnarray*}
(\bar{x}_n)^n &=& |f(\bar{x}_n)| = g(\bar{x}_n) - g(0) \le \, \int_0^{{\!\!\!\!\!\!\!*} \, {\,\,\,\,\bar{x}_n}} g'(x) \, dx \le \,  \int_0^{{\!\!\!\!\!\!\!*} \, {\,\,\,\,\bar{x}_n}} |g'(x)| \, dx \\
&=&  \int_0^{{\!\!\!\!\!\!\!*} \, {\,\,\,\,\bar{x}_n}} |f'(x)| \, dx \le \, \int_0^{\bar{x}_n} C \, \frac{|f(x)|}{x} \, dx < C \int_0^{\bar{x}_n} \frac{x^n}{x} \, dx = \frac{C}{n} (\bar{x}_n)^n < (\bar{x}_n)^n.
\end{eqnarray*}
This tells us that $E_n = \emptyset$ i.e. 
\[
\forall n > C, \qquad  |f(x)| < x^n \qquad \forall x \in ]0,1].
\]
This entails that $|f(x)| \equiv 0.$ 
\end{proof}

\end{document}